# Minimum Time Escape from a Circular Region of a Dubins Car

Isaac E. Weintraub
*Control Science Center*
*Air Force Research Laboratory*
Wright-Patterson AFB, OH 45433
isaac.weintraub.1@afrl.af.mil

Alexander Von Moll
*Control Science Center*
*Air Force Research Laboratory*
Wright-Patterson AFB, OH 45433
alexander.von_moll@afrl.af.mil

Meir N. Pachter
*Electrical and Computer Engineering*
*Air Force Institute of Technology*
Wright-Patterson AFB, OH 45433
meir.pachter@afit.edu

*Abstract—* A turn-constrained evader strives to escape a circular region in minimum time.

*Index terms*—Path Planning, Optimization, Robotics, Navigation

## I. Introduction

The importance of path-planning and vehicle routing has been an important area of research for centuries. In optimal path planning, the objective is for a vehicle to obtain a policy to accomplish an objective/criteria as best as possible. In this paper, the objective is for a mobile turn-constrained vehicle, called an evader, to escape a circular domain in minimum time. The initial location, pose, and turn constraints of the evader dictate the strategy used to escape the circular region.

The time-optimal motion of a turn-constrained vehicle between two points was demonstrated in [1], [2]; optimal trajectories were shown to consist of combinations of arc-turns of minimum radius and straight lines. Dubins' vehicle model and results have been leveraged by countless others for a broad host of applications. Later, extensions were made to time-optimal controls for turn-constrained vehicles that can go forward or backward [3] and generalizations to turn-constrained vehicles [4].

Applications for this work extend to surveillance evasion such as work in [5]–[8]. Also related are scenarios of border patrol of turn-constrained vehicles such as [9]. The shortest path for a vehicle to reach a circular region tangentially from the outside was investigated in [10]–[12].

In this paper, the use of the method of characteristics is leveraged [13]. This is akin to solving the optimal control problem, working backward from a terminal manifold related to optimal control [14], dynamic programming [15], and Pontryagin's minimum principle [16]. A related paper that leverages the method of characteristics for solving optimal control is found in [17].

Time-optimal strategies for a Dubins vehicle to navigate through a set number of points have also been considered. In [18] strategies for navigating through three points are provided. More generalized approaches for more than three points have been considered in [19], [20], to visit a set of regions [21], in cliques [22], and fluid currents [23], [24].

In a very related work, a Dubins vehicle strives to escape a circular region in minimum time [25]. This paper considers the very same problem; but, reduces the state-space size by considering the reduced space of dimension 2 rather than the Cartesian coordinates which require 3. Moreover, we report the time it takes to escape and the conditions that dictate which strategy to take as a function of initial conditions. The analysis used herein also leverages the method of characteristics to argue for the usage of geometry.

## II. Problem Definition

In this classical problem, the objective is for an evader, whose position is defined as $\mathbf{x} = (r, \theta)^\top \in \mathbb{R}^2$ to escape the circular region $\mathcal{C}$ centered at $C = (0,0)$ and has radius $\rho$, in minimum time. Let $r$ be distance of the evader from $C$, and defined over the domain $0 \leq r \leq \rho$, and let $\theta$ be the heading angle of the evader with respect to the $\hat{r}$-axis which points outward from $C$, and is defined over $-\pi \leq \theta \leq \pi$. The objective is the min-time escape:

$$\operatorname*{argmin}_{u \in \mathcal{U}}\{J\} = \operatorname*{argmin}_{u \in \mathcal{U}}\left\{\int_0^{t_f} 1 \, \mathrm{d}t\right\} = \operatorname*{argmin}_{u \in \mathcal{U}}\{t_f\} \quad (1)$$

where the final time $t_f$ is defined as

$$t_f = \{t \mid r(t) = 1, -\pi \leq \theta(t) \leq \pi\} \quad (2)$$

Rather than consider a non-holonomic vehicle, that can turn in-place, the evader, $E$, is modeled having a minimum turn radius, $R$, and moves with speed $v_E$. A figure that graphically describes $E$ and its escape from the circular region $\mathcal{C}$ is depicted in Figure 1.

The state space dynamics, which describe the motion of $E$, is represented by a set of first-order ordinary differential equations: $\dot{\mathbf{x}} = \boldsymbol{f}(\mathbf{x}(t), u(t), t)$,

---
This work was funded in part by LRIR 24RQCOR002.
Distribution Statement A. Approved for public release: distribution is unlimited. AFRL-2024-2504, CLEARED 08 MAY 2024

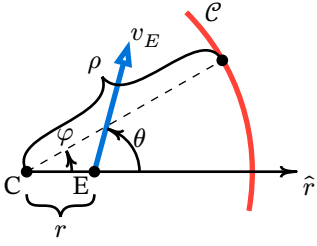

Figure 1: An evader, $E$, strives to escape a circular region, $\mathcal{C}$, in minimum time.

$$\dot{r} = v_E \cos\theta, \quad r(0) = r_0(>0)$$
$$\dot{\theta} = \frac{v_E}{R}u - \frac{v_E}{r}\sin\theta, \quad \theta(0) = \theta_0, \quad 0 \le t. \quad (3)$$

It is assumed that $E$ starts some initial notional distance from the center of the circular region $r_0 > 0$ and has an initial heading $\theta_0$. Without-loss of generality, it is assumed that the initial heading of the evader $0 \le \theta_0 \le \pi$. By symmetry, the approaches taken and presented in this work can be reflected about the $\hat{r}$-axis. For this reason, this paper focuses attention to the top-half plane.

The evader operates in the real plane and a feedback strategy for the evader is desired as a function of the two states: $\mathbf{x} = (r, \theta)^\top \in \mathbb{R}^2$. By scaling the problem using the radius of the circular region, $\rho$, the problem can be further generalized as follows:

$$r \to \tfrac{r}{\rho}, \quad R \to \tfrac{R}{\rho}, \quad \rho \to \tfrac{\rho}{\rho} = 1 \quad (4)$$

The generalized equations of motion (dynamics) from (3) when scaled by (4) become

$$\dot{r} = \cos(\theta), \quad r(0) = r_0(>0)$$
$$\dot{\theta} = -\tfrac{1}{r}\sin\theta + \tfrac{1}{R}u, \quad \theta(0) = \theta_0, \quad 0 \le t. \quad (5)$$

The terminal manifold $\mathcal{M} = \{(r,\theta) \mid r = 1, -\pi \le \theta \le \pi\}$. That is, the terminal manifold is the circular boundary where the evader strives to reach in minimum time. The two-dimensional state space $(r,\theta)$ is a cylindrical manifold.

## III. METHODS AND ANALYSIS

The approach for solving this problem is two-part: first, the method of characteristics is leveraged to demonstrate the monotonic nature of critical variables. This fact is then used to justify the use of geometric methods for obtaining the optimal solution.

### A. Method of Characteristics

By the method of characteristics, the Usable Part (UP) of the terminal manifold $\mathcal{M}$, can be obtained by taking the inner-product of the state dynamics (5) with the inward pointing normals to $\mathcal{M}$. The inward pointing normals are $\vec{n} = (-1, 0)^\top$. The UP is defined as

$$\mathrm{UP} \triangleq \max_{u(t) \in \mathcal{U}} \boldsymbol{f}(\mathbf{x}(t), u(t), t) \cdot \vec{n} < 0, \quad (6)$$

where $\dot{\mathbf{x}} = \boldsymbol{f}(\mathbf{x}(t), u(t), t)$ is a vector of the dynamics as defined in (5) and this inner-product is maximized for the admissible control $u(t) \in \mathcal{U}$. The UP is therefore

$$\mathrm{UP} \triangleq \max_{-1 \le u(t) \le 1} (-1 \cdot \cos\theta) < 0 \quad (7)$$

By (7) the UP is the right-half of the circular boundary:

$$\mathrm{UP}(\mathcal{M}) = \left\{(r, \theta) \mid r = 1, -\tfrac{\pi}{2} < \theta < \tfrac{\pi}{2}\right\} \quad (8)$$

In this problem, there exists a Universal Line (UL), which is the line segment described as

$$\mathrm{UL} \triangleq \{(r, \theta) \mid 0 < r < 1, \theta = 1\} \quad (9)$$

On the UL, as shown in (9), the costate, $\lambda_\theta = 0$; by symmetry. The UL represents a singular surface (or line in this case) by which tributary trajectories reach prior to reaching the terminal manifold.

The Hamiltonian for this problem is:

$$\mathcal{H} = -1 + \lambda_r \cos\theta + \tfrac{1}{R}\lambda_\theta u - \tfrac{1}{r}\lambda_\theta \sin\theta \quad (10)$$

From Pontryagin's minimum:

$$\mathcal{H}(\mathbf{x}^*(t), u^*(t), \lambda^*(t), t) \le \mathcal{H}(\mathbf{x}^*(t), u(t), \lambda^*(t), t) \quad (11)$$

From (11), the following is attained:

$$u^*(t)\lambda_\theta \le u(t)\lambda_\theta \quad (12)$$

By (12), the optimal control is:

$$u^* = \mathrm{sign}(\lambda_\theta), \quad (13)$$

where

$$\mathrm{sign}(x) = \begin{cases} 1 & x > 0 \\ 0 & x = 0 \\ -1 & x < 0 \end{cases} \quad (14)$$

The costate dynamics are attained using the partials of the Hamiltonian with respect to the state variables:

$$\dot{\lambda}_r = -\tfrac{\partial \mathcal{H}}{\partial r} = -\tfrac{\lambda_\theta \sin\theta}{r^2},$$
$$\dot{\lambda}_\theta = -\tfrac{\partial \mathcal{H}}{\partial \theta} = \lambda_r \sin\theta + \tfrac{\lambda_\theta \cos\theta}{r} \quad (15)$$

The Euler-Lagrange equations (dynamics) in retrograde form are:

$$\mathring{r} = -\cos\theta \qquad r(\tau = 0) = 1$$
$$\mathring{\theta} = \tfrac{1}{r}\sin\theta - \tfrac{1}{R}\mathrm{sign}(\lambda) \qquad \theta(\tau = 0) = t_f$$
$$\mathring{\lambda}_r = \tfrac{1}{r^2}\lambda_\theta \sin\theta \qquad \lambda_r(\tau = 0) = \sec\theta_f$$
$$\mathring{\lambda}_\theta = -\lambda_r \sin\theta - \tfrac{1}{r}\lambda_\theta \cos\theta \qquad \lambda_\theta(\tau = 0) = 0, 0 \le \tau$$
$$(16)$$

Substitution of (13) into the Hamiltonian in (10), the optimal Hamiltonian is

$$\mathcal{H}^* = -1\lambda_r \cos\theta + \tfrac{1}{R}|\lambda_\theta| - \tfrac{1}{r}\lambda_\theta \sin\theta \equiv 0 \quad (17)$$

Therefore on the UL, $\lambda_r = 1$. Furthermore, from construction, $\lambda_\theta = -V_\theta$, the costate value of $\theta$ is negative of the sensitivity of the value function (time to escape) with respect to the variable $\theta$. One may observe that for smaller values of $\theta$ are implicitly related to faster escape times. A proof of this is omitted due to the limited space available in this paper and is left for future work.

The property of $|\theta|$ monotonically decreases until it converges to zero; whereby $\mathbf{x}(t)$ reaches the UL or the terminal manifold (whichever occurs first) and gives rise to the optimal strategy – motivating the use of geometry: 1) Turn until reaching the UL, whereby the trajectory is a straight-line or 2) Escape the circular region (by way of turning the whole time) before reaching the UL.

### B. Geometric Analysis

The evader is located at $E = (r \cdot \hat{x}, 0 \cdot \hat{y})$ at time $t_0$ and has a heading $\theta$ as described in Figure 1. Because the optimal strategy of $E$ is to execute a turn until it reaches the UL, after which, it holds course and escapes the circular region; the use of geometry (circles, lines, and tangents), can be used to describe the time-optimal strategy for $E$. As previously shown, three strategies exist for $E$: straight, turn-straight, and turn.

$$u_E = \begin{cases} 0 & \theta = 0 \\ -1 & \theta > 0 \end{cases} \quad (18)$$

1) *Straight Line:* When $\theta(t_0) = 0$, the evader is pointed directly toward the circular boundary and lies on the UL. In this case, the optimal strategy is for $E$ to hold course until it exits the circular boundary. This condition is pictorially described in Figure 2.

The time for the evader to escape the circular region is

$$t_f = t_{\text{straight}} = \frac{\rho - r}{v_E} \quad (19)$$

2) *Turn-Straight:* The turn-straight condition occurs when a tangent point $T$ is inside the circular escape region. If one defines $T$ in Cartesian coordinates as $T = (x_T, y_T)$, then $T$ is inside the circular region when $x_T^2 + y_T^2 < \rho^2$. This condition is pictorially described in Figure 3.

For this portion, the top-half plane is considered and therefore only right-turns are considered. The turn-circle of the evader is centered at a point, $O$,

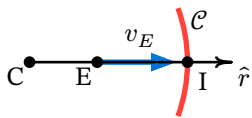

Figure 2: The evader, $E$, escapes a circular region, $\mathcal{C}$ in minimum time and lies on UL at the onset: $\theta(t_0) = \theta(t_f) = 0$.

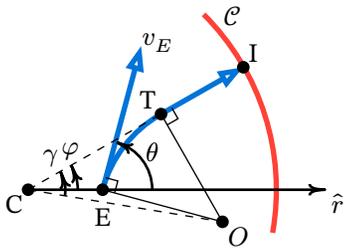

Figure 3: The evader, $E$, escapes a circular region, $\mathcal{C}$, in minimum time and reaches the UL prior to escape, $\theta(t_f) = 0$.

$$O = (r\hat{x} + R\sin(\theta)\hat{x}, -R\cos(\theta)\hat{y})^\top. \quad (20)$$

The current objective is to locate $T$ so that one can determine if $E$ takes a straight-turn strategy. Define the angle $\gamma$ as shown in Figure 3,

$$\varphi = \sin^{-1}\left(\frac{R}{\|O\|}\right) - \sin^{-1}\left(\frac{R\cos\theta}{\|O\|}\right) \quad (21)$$

Consider the vector sum:

$$\vec{CT} + \vec{TO} = \vec{CE} + \vec{EO} \Rightarrow$$
$$\|T\|\cos\varphi\hat{x} + \|T\|\sin\varphi\hat{y} + R\sin\varphi\hat{x} - R\cos\varphi\hat{y} \quad (22)$$
$$= r\hat{x} + R\sin\theta\hat{x} - R\cos(\theta)\hat{y}$$

Separating (22), the pair of equations can be obtained:

$$\|T\|\cos\varphi + R\sin\varphi = r + R\sin\theta \quad (23)$$
$$\|T\|\sin\varphi - R\cos\varphi = -R\cos\theta \quad (24)$$

Substitution of $\varphi$ from (21) into (24),

$$\|T\|\sin(\sin^{-1}(a) - \sin^{-1}(b)) \\ -R\cos(\sin^{-1}(a) - \sin^{-1}(b)) = -R\cos(\theta) \quad (25)$$

where the intermittent variables are

$$a = \frac{R}{\|O\|}, b = \frac{R\cos\theta}{\|O\|} \quad (26)$$

further expansion of (25) and solving for $\|T\|$:

$$\|T\| = R\frac{\cos(\sin^{-1}(a) - \sin^{-1}(b)) - \cos\theta}{\sin(\sin^{-1}(a) - \sin^{-1}(b))} \quad (27)$$

substitution of the intermittent variable in (26) into (27), the following is obtained:

$$\|T\| = \frac{\sigma_1\sigma_2 - \sigma_2^2\cos\theta}{\sigma_1 - \sigma_2} \quad (28)$$

where

$$\sigma_1 = \sqrt{r^2 + R^2 + 2rR\sin\theta - R^2\cos^2\theta}, \\ \sigma_2 = \sqrt{r^2 + 2rR\sin\theta} \quad (29)$$

The tangent point, $T$, is inside the the circular region $\mathcal{C}$ if $\rho - \|T\| > 0$. In this case, the time it takes to escape the circular region is:

$$t_f = t_{\text{turn-straight}} = \underbrace{\frac{R(\theta - \varphi)}{v_E}}_{\text{turn}} + \underbrace{\frac{\rho - \|T\|}{v_E}}_{\text{straight}} \quad (30)$$

3) *Turn-Only:* The turn-only condition occurs when $T$ is outside or on $\mathcal{C}$. Let $T = (x_T, y_T)$, then $T$ is outside or on $\mathcal{C}$ when $x_T^2 + y_T^2 \geq \rho^2$. This condition is described in Figure 4.

The intercept point, $I$, where the evader escapes $\mathcal{C}$ is dictated by the circle-circle intercept between the turn-circle and $\mathcal{C}$. The intercept point is located (in polar coordinates) at

$$I \triangleq (\rho\hat{r}, \alpha\hat{\theta}) \quad (31)$$

Using the law of cosines over the $\triangle CIO$:

$$R^2 = \|O\|^2 + \rho^2 - 2\|O\|\rho\cos\delta \quad (32)$$

Solving (32) for the angle $\delta$,

$$\cos\delta = \frac{\|O\|^2 + \rho^2 - R^2}{2\rho\|O\|} \quad (33)$$

By construction the angle $\alpha$ that locates $I$ is

$$\alpha = \cos^{-1}\left(\frac{\|O\|^2 + \rho^2 - R^2}{2\rho \|O\|}\right) - \sin^{-1}\left(\frac{R\cos\theta}{\|O\|}\right) \quad (34)$$

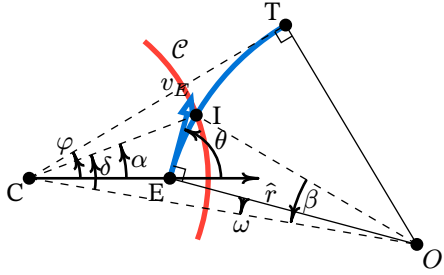

Figure 4: The evader, $E$, escapes a circular region, $\mathcal{C}$, in minimum time and, due to the large turn-radius, is unable to reach the UL prior to escape, $\theta(t_f) > 0$.

Therefore the intercept point may be written as:
$$I = (\rho\cos\alpha\,\hat{x},\, \rho\sin\alpha\,\hat{y}) \quad (35)$$

Consider $\triangle IOC$; the angle $\beta \triangleq \angle IOC$, using the law of cosines:
$$\rho^2 = \|O\|^2 + R^2 - 2\|O\|R\cos(\beta) \quad (36)$$

Solving for $\cos(\beta)$
$$\cos(\beta) = \frac{\|O\|^2 + R^2 - \rho^2}{2\|O\|R} \quad (37)$$

Also, consider $\triangle EOC$; the angle $\omega \triangleq \angle EOC$, using the law of cosines again,
$$r^2 = \|O\|^2 + R^2 - 2\|O\|R\cos(\omega) \quad (38)$$

solving for $\cos(\omega)$
$$\cos\omega = \frac{r^2 - \|O\|^2 - R^2}{2\|O\|R} \quad (39)$$

The angle traversed by $E$ from its initial condition until it reaches $I$ (i.e. the turn circle arc that is contained inside the circular region) is obtained by subtracting $\beta - \omega$. Therefore the time taken to leave the circular region is
$$t_{\text{turn}} = \frac{R}{v_E}(\beta - \omega) \quad (40)$$

More explicitly, (40) can be written as
$$t_f = \frac{R}{v_E}\left[\cos^{-1}\left(\frac{\|O\|^2 + R^2 - \rho^2}{2\|O\|R}\right) - \cos^{-1}\left(\frac{r^2 - \|O\|^2 - R^2}{2\|O\|R}\right)\right] \quad (41)$$

## IV. Conclusions

In this paper, the minimum-time escape of a turn-constrained evader was considered. The time-optimal strategies for the evader were obtained using geometry, motivated by the method of characteristics which stipulated that optimal strategies consisted of straight, turn-straight, and turn. Using geometry, the time until escape was reported for each of the three possible strategies. Future work could consider multi-agent scenarios incorporating this strategy, non-circular regions of space, and more attention to nuances when leveraging the method of characteristics.